\documentclass[a4paper]{article}

\usepackage{amsfonts}
\usepackage{amssymb}
\usepackage{amsmath} 
\usepackage{amsthm}
\usepackage{tikz}
\usepackage{color}
\usepackage[ruled, lined, longend, linesnumbered]{algorithm2e}
\usepackage{authblk}

\definecolor{grey}{gray}{0.92}

\newcommand*{\N}{\ensuremath{\mathbb{N}}}

\theoremstyle{nummermitklammern}
\newtheorem{defsatzusw}{}[section]
\newtheorem{definition}[defsatzusw]{Definition}
\newtheorem{theorem}[defsatzusw]{Theorem}

\newtheorem{conjecture}[defsatzusw]{Conjecture}
\newcounter{nummer}
\begin{document}

\author[1]{Robert Scheidweiler
\thanks{scheidweiler@math2.rwth-aachen.de}}
\affil[1]{\small{Paderborn, Germany}}
\author[2]{Eberhard Triesch
\thanks{triesch@math2.rwth-aachen.de}}
\affil[2]{Lehrstuhl 2 f\"ur Mathematik, RWTH Aachen University, Aachen, Germany}

\title{Upper and Lower Bounds for Competitive Group Testing}

\date{}

\maketitle

\begin{abstract}
 We consider competitive algorithms for adaptive group testing problems. In the first part of the paper, we develop an algorithm with competitive constant $c<1.452$ thus improving the up to now best known algorithms with constants $1.5+\epsilon$ from 2003. In the second part, we prove the first nontrivial lower bound for competitive constants, namely that $c$ is always larger than $1.31$. \\
\end{abstract}

\section{Introduction}
\noindent

Assume that we are given some finite set $S$, $|S|=n$, of {\it items} which can be
{\it good} or {\it defective}. Our task is to identify the set $D$ of defectives by
(successively) choosing subsets $W_1,W_2,\ldots $ of $S$ and testing whether  $D$ is disjoint from $W_i$ or not. In choosing $W_i$, we may assume that the results of the first $i-1$ tests are already known ({\it adaptive} search). We want to
minimize the number of tests needed to identify the set $D$ in the worst case.

This problem is called the {\it adaptive group testing} problem and dates back to 1943 when Dorfman \cite{Dorfman1943} in his seminal paper discussed possibilities to
reduce the number of blood tests for syphilitic antigene of American inductees.
Although group testing was not used in this context, many applications in different areas such as quality control, information theory, computer networks and clone library were found in the following years. The most complete account on
group testing can be found in \cite{Du2000}. %(Du Hwang 2000)
The books \cite{Aigner1988} % (Aigner) 
and
\cite{Ahlswede1979} % (Ahlswede Wegener) 
treat search problems in more generality but also cover the subject while
\cite{Du2006} %(Du Hwang 2006)
 focusses on nonadaptive group testing with applications in the life sciences. In this paper, we deal with adaptive group
testing only.

\noindent
Most papers on this subject considered the case that the number $d$ of defectives is known in advance, the so-called $(d,n)$-group testing problem.
We write $M_A(d,n)$ for the maximum number of tests required by a group testing algorithm $A$ to identify all items and denote by $M(d,n):= \min_{A} M_A(d,n)$ the worst case number of tests for the problem.
By general
considerations, we have the {\it information theoretic lower bound}

$$M(d,n) \geq \left\lceil \log {n \choose d} \right\rceil \geq d \log \frac{n}{d},$$

\noindent
where the base of $\log$ is two throughout this paper.

\noindent
There is a famous and beautiful conjecture on the so-called cutoff point for adaptive group testing due to Hu, Hwang and Wang:
If the ratio $r:=n/d$ is at most $3$, it is optimal to test each item individually.

\begin{conjecture}\cite{Hu1981}\label{Corollary:Cutoff}
For $0 < d < n$ and $r \leq 3$, 
\begin{align*}
M(d,n)=n-1.
\end{align*} 
\end{conjecture}
\noindent

The equation $M(d,n)=n-1$ was proved if $r\leq \frac{21}{8}=2.625$ by Du and Hwang  \cite{Du1982}.
\noindent
Leu, Lin and Weng extended this to $r\leq \frac {43}{16}=2.6875$ for $d\geq 193$ \cite{Leu2002}
\noindent
and Riccio and Colbourn to $r< \log_{\frac{3}{2}}3\approx 2.709$ for sufficiently large $d$
depending on $r$ \cite{Riccio2000}.
\noindent
Wegener, Fischer and Klasner proved the conjecture if the cardinality of the test sets is
at most 2 \cite{Fischer1999}. For later use, we record the result of Leu, Lin and Weng
\begin{theorem}\label{Leu2002}\cite{Leu2002}
For $\frac{16}{43}n\leq d < n$, $d \geq 193$,
\begin{align*}
 M(d,n)=n-1.
\end{align*} 
\end{theorem}

\noindent
For $r>\frac{43}{16}$, the best group testing algorithms we know about are due to A. Allemann \cite{Allemann2013}:

\begin{theorem}\label{Allemann2013}\cite{Allemann2013}
There is an algorithm which, for all $0<d\leq \frac{n}{2}$,
finds $d$ defectives by at most 
$$\log{n\choose d}+0.255d+\frac{1}{2}\log d +6.5$$
tests.
\end{theorem}

\noindent In practice, the prior knowledge of $d$ is often an unrealistic assumption.
However, if nothing is known about $d$, the worst case number of tests is clearly $n$ and can be realized by testing all elements
individually. Suppose the number of defectives turns out to be one. Then those $n$ singleton tests are very many as compared to the $\lceil\log n\rceil$
tests needed if $d=1$ is known in advance. Du and Hwang thus asked for an algorithm $A$ with the following property: 
Algorithm $A$ successfully identifies $D$ without prior knowledge of $d$ for all $n$. Denote by $M_A(d\, |\, n)$ the worst case number of tests used by $A$ if
the number of defectives turns out to be $d$. Then the ratio $M_A(d\, |\, n)/M(d,n)$ should be bounded by some constant c.
More precisely:

\noindent
\begin{definition}\cite{Du1993}
Let $c$ be some real number, $c\geq 1$. Algorithm $A$ is called $c$-competitive if there is some constant $a \in \N$ such that 
\begin{align*}
 M_A(d\,|\,n) \leq c M(d,n)+a \ ,
\end{align*}
for all $0 \leq d < n$. The number $c$ is called a competitive ratio for $A$.
\end{definition}

\noindent
Du and Hwang gave the first competitive algorithm with $c= 2.75$ in 1993 \cite{Du1993}. Several improvements were given by Bar-Noy, 
Hwang, Kessler and Kutten $(c=2)$ in 1994 \cite{Bar-Noy1994}, by Du, Xue, Sun and Cheng $(c=1.65)$ in 1994 \cite{Du1994} and by Schlaghoff
and Triesch $(c=1.5+\epsilon)$ in 2005 \cite{Schlaghoff2003}.\\
In  \cite[p.~77]{Du2000} %p.77,
  we read the following remark: 
\textit{``However, it certainly requires a new technique in order to push the competitive ratio down under $1.5$.''}

No nontrivial \textit{lower bound} on the competitive ratio has been proved so far. For an interesting related paper, see \cite{Damaschke2009}.

In the next section, we are going to present an algorithm with competitive ratio $c< 1.452$. We then go on to prove that for each $c$-competitive algorithm we have $c>1.31$.

\section{A 1.452-competitive algorithm}
In this section, we present a new algorithm $A$ with competitive ratio $c< 1.452$.
Before we start explaining the ideas, we note two useful inequalities.
The first is an estimate of  $\log \left({n \choose d}\right)$ by Stirling's formula
 $k!= \sqrt{2\pi k} \left(\frac{k}{e}\right)^k \theta_k$ with $1 < \theta_k \leq e^{\frac{1}{12k}}$.
Recall that $r$ denotes the quotient $n/d$.

\begin{theorem}\cite{Allemann2003}\label{Theorem: Stirling}
For $0 < d \leq \frac{n}{2}$ we have:
\begin{align*}
M(d,n) \geq d \cdot \left( \log r+(r-1)\log \left( \frac{r}{r-1}\right)\right)-\frac{1}{2} \log d -\frac{3}{2} \ .
\end{align*} 
For $r \geq 3$ we have $(r-1)\log \left( \frac{r}{r-1}\right)> 1.1699$.
\end{theorem}
\noindent 
For a detailed proof, we refer the reader to \cite{Allemann2003}.\\

\noindent The second tool we need is the following standard inequality which is implied by the \textit{concavity of the logarithm} and \textit{Jensen's inequality}:
\begin{theorem}\label{Jensen}
Let $0 < d_i \leq n_i$ for $1 \leq i \leq k$, $d:= \sum\limits_{i=1}^k d_i$, $n:= \sum\limits_{i=1}^k n_i$. Then the following inequality holds:
\begin{align*}
 \sum\limits_{i=1}^k d_i\log \frac{n_i}{d_i} \leq d \log \frac{n}{d}.
\end{align*}
\end{theorem}

\noindent
The algorithm we present naturally decomposes into ``phases'' where some good and some defective items are identified. If in phase $i$ we detect $g_i$ good and $d_i$ defective items, we will show that the corresponding
number of tests $t_i$ is at most $ 1.45198 d_i (\log \frac{g_i+d_i}{d_i}+1.1699 )$. Summing over $i$ and estimating by Theorem \ref{Jensen}
we obtain an upper bound of $ 1.45198\,  d(\log \frac{n}{d}+1.1699 )$ for the total number of tests.

%%%%%%%%%%%%%%%%%%%%%%%%%%%%%%%%%%%%%%%%%%%%%%%%%%%

\noindent
For a formal description compare Algorithm \ref{Main} which uses some subalgorithms described below. The variable $S$ denotes the set of still unclassified elements through the execution of the algorithm. The command TEST($X$) means that we test whether $X$ contains at least one defective element or not. In the first case, the test is called
positive and negative in the second. 

Procedure DIG (Algorithm \ref{Dig}) is applied to a defective set $X$ and finds a defective element by binary splitting. The number of tests is of course $\lceil \log |X| \rceil$.

More interesting is Algorithm \ref{Fourtest} called FOURTEST. It is also applied to some defective set $X$ where we assume that $|X|=2^l$ is a power of two, $|X|\geq 4$. It first splits $X$ into four disjoint subsets $T,U,V,W$ of equal size and tests one of them, say $T$. If it is defective, DIG is applied to find a defective element by $\log |T|=l-2$ additional tests. Otherwise, the sets $U$ and $V$ are tested. If one or both of them are defective, DIG is applied to find one or two defective elements. If none of them is defective, we know that $W$ is defective and, again find a defective by applying DIG. Summarizing, FOURTEST produces one the following results:
\begin{itemize}
\item one defective is identified using $l-1$ tests.
\item one defective and $2^{l-1}$ good elements are detected by $3+l-2=l+1$ tests.
\item two defectives and $2^{l-2}$ good elements are detected by $3+2(l-2)=2l-1$ tests.
\item one defective and $3\cdot 2^{l-2}$ good elements are found by $3+l-2=l+1$ tests.
\end{itemize}

%The subalgorithms $A_k$ with $2\leq k\leq 6$ are explained below.

We now begin discussing our main procedure, i.e., Algorithm \ref{Main}:

\subsubsection*{Main Procedure}
At the beginning, we check whether the number of elements in $S$ is larger than $576=8+8\cdot 64$. If less elements are unclassified, we use the command INDVTEST($S$) which stands for individual tests of all elements of $S.$ The constant $a$ in the definition of competitiveness takes care of this constant number of tests. (A similar remark applies to the subalgorithms $A_k$, $2\leq k\leq 6,$ compare Algorithm \ref{Ak}.) If enough elements are unclassified, we test a subset of cardinality $64$ and test it. Let us suppose first that this group of 64 elements is good.

%\begin{algorithm}
%	\caption{MAIN (set $S$) }
%	\label{Main}
%	\textbf{global } $S$\\
%	\textbf{global } $D = \emptyset$\\
%	\textbf{global } $G = \emptyset$ \\
%	$i=0$ \\
%	\While{$|S|> 576$}{
%		$X = \min(64\cdot 8^i, |S|)$ elements from $S$ \\
%		\eIf{$|X|=64\cdot 8^i$}{
%		TEST$(X)$ \\
%		\eIf{$X$ is good}{%
%			$G = G \cup X$ \\
%			$S = S \setminus X$ \\
%			$i = i+1$  \\
%			\If{$i = 2$}{%
%				TEST$(S)$\\
%				\If{$S$ is good}{
%					$G = G \cup S$ \\
%					$S = \emptyset$ \\
%				}
%		}}{
%				\eIf{$i=0$}{ 
%					$A_6(X,S\setminus X)$}
%				{
%					FOURTEST$(X)$  \\
%					$i=0$  }
%	}}{DIG$(X)$}}
%	INDVTEST$(S)$\\
%	\textbf{return} $D,$ $G$
%\end{algorithm}
%

\begin{algorithm}
	\caption{MAIN (set $S$) }
	\label{Main}
	\textbf{global } $S$\\
	\textbf{global } $D = \emptyset$\\
	\textbf{global } $G = \emptyset$ \\
	$k=0$ \\
	\While{$|S|> 576$}{
		$X = \min(64\cdot 8^k, |S|)$ elements from $S$ \\
		\eIf{$|X|=64\cdot 8^k$}{
			TEST$(X)$ \\
			\eIf{$X$ is good}{%
				$G = G \cup X$ \\
				$S = S \setminus X$ \\
				\If{$k = 1$}{%
					TEST$(S)$\\
					\If{$S$ is good}{
						$G = G \cup S$ \\
						$S = \emptyset$ \\
					}

			}$k= k+1$}{
				\eIf{$k=0$}{ 
					$A_6(X,S\setminus X)$}
				{
					FOURTEST$(X)$  \\
					$k=0$  }
		}}{DIG$(X)$}}
	INDVTEST$(S)$\\
	\textbf{return} $D,$ $G$
\end{algorithm}

\begin{algorithm}
	\caption{FOURTEST (set $X$) }
	\label{Fourtest}
	$T,U,V,W$ = a partition of $X$ into $4$ sets of equal cardinality \\
	TEST($T$)\\
	\eIf{$T$ is defective} 
	{${DIG}(T)$\\}
	{ 
	$G = G \cup T$ \\
	$S = S\setminus T$\\
	TEST($U$) \\
	TEST($V$) \\
	\eIf{$U$ is defective}{
	${DIG}(U)$} 
	{$G = G \cup U$ \\
	$S = S\setminus U$}
	\eIf{$V$ is defective}{
	${DIG}(V)$} 
	{$G = G \cup V$ \\
	$S = S\setminus V$} 
	\If{$U$ and $V$ are good}{
	DIG($W$)}}
\end{algorithm}

%\begin{figure}
%\begin{center}
%{\footnotesize\colorbox{grey}{\parbox{11cm}{
%{\normalsize\textbf{Procedure FOURTEST(set $X$)}} \\[0.25cm]
%$T,U,V,W$ = a partition of $X$ into $4$ sets of equal cardinality \\
%TEST($T$)\\[0.25cm]
%\textbf{if} $T$ is defective:\\ 
%\tab${DIG}(T)$ \\
%\textbf{else}: \\ 
%\tab $G = G \cup T$ \\
%\tab $S = S\setminus T$\\[0.25cm]
%\tab TEST($U$) \\
%\tab TEST($V$) \\[0.25cm]
%\tab \textbf{if} $U$ is defective: \\
%\tab[2cm] ${DIG}(U)$\\
%\tab \textbf{else}: \\ 
%\tab[2cm] $G = G \cup U$ \\
%\tab[2cm] $S = S\setminus U$\\
%\tab \textbf{if} $V$ is defective: \\
%\tab[2cm] ${DIG}(V)$\\
%\tab \textbf{else}: \\ 
%\tab[2cm] $G = G \cup U$ \\
%\tab[2cm] $S = S\setminus U$\\
%\tab \textbf{if} $U$ and $V$ are good: \\
%\tab[2cm] DIG($W$) \\[0.25cm]
%
%}}}
%\caption{Procedure \textbf{FOURTEST}}
%\label{Fourtest}
%\end{center} 
%\end{figure}

\begin{algorithm}
	\caption{DIG (set $X$) (compare \cite{Du2000}) }
	\label{Dig}

	\While{$|X|\geq 1$}{
	$Y = \lceil |X| / 2  \rceil$ elements from $X$\\
	TEST($Y$)\\
	\eIf{$Y$ is good}{
			$G = G \cup Y$\\
			$S = S \setminus Y$\\
			$X = X \setminus Y$\\} 
			{$X = Y$}
			
		 }
	$S = S \setminus X$\\
	$D = D \cup X$			
	
\end{algorithm}

%\begin{figure}
%\begin{center}
%{\footnotesize\colorbox{grey}{\parbox{8.5cm}{
%{\normalsize\textbf{Function DIG(set $X$)}} \\[0.25cm]
%\textbf{while} $|X|\geq 1$:\\
%\tab $Y = \lceil |X| / 2  \rceil$ elements from $X$\\
%\tab TEST($Y$)\\
%\tab \textbf{if} $Y$ is good:\\
%\tab[2cm] $G = G \cup Y$\\
%\tab[2cm] $S = S \setminus Y$\\
%\tab[2cm] $X = X \setminus Y$\\ 
%\tab \textbf{else}:\\
%\tab[2cm] $X = Y$ \\[0.25cm]
%$S = S \setminus X$\\
%$D = D \cup X$
%
%}}}
%\caption{Function \textbf{DIG} (see \cite{Du2000})}
%\label{Dig}
%\end{center} 
%\end{figure}

\begin{algorithm}
	\caption{$\text{A}_k (\text{set } X_k, \text{set } T_k  )$ for  $3\leq k\leq 6$ }
	\label{Ak}
	\textbf{global } $X_k$\\
	\textbf{global } $T_k$\\
	\eIf{$|T_k|<2^{k-1}$}{
	INDVTEST$(T_k)$}
	{$Y = 2^{k-1}$ items from $T_k$\\
	TEST($Y$)\\
	\eIf{$Y$ is good}{ 
	$G = G \cup Y$ \\
	$S = S \setminus Y$ \\
	FOURTEST($X_k$) \\
	\If{$3\leq k \leq 5$}
	{$\text{A}_{k+1}(X_{k+1},T_k\setminus Y)$\\}}
	{\eIf{$k>3$}
	{$\text{A}_{k-1}(Y,T_k\setminus Y)$ \\}
	{
	 $\text{A}_2(Y)$\\}}}

\end{algorithm}

%\begin{figure}
%\begin{center}
%{\footnotesize\colorbox{grey}{\parbox{8cm}{
%{\normalsize \textbf{Procedure $\text{A}_k (\text{set } X_k, \text{set } T_k  )$} \tab $3\leq k\leq 6$}\\[0.25cm]
%
%\textbf{global variable} $X_k$\\
%\textbf{global variable} $T_k$\\[0.25cm]
%
%\textbf{if} $|T_k|<2^{k-1}:$ \\
%\tab INDVTEST$(T_k)$\\
%\textbf{else}:\\
%\tab $Y = 2^{k-1}$ items from $T_k$\\
%\tab TEST($Y$)\\
%\tab \textbf{if} $Y$ is good:\\ 
%\tab[2cm] $G = G \cup Y$ \\
%\tab[2cm] $S = S \setminus Y$ \\
%\tab[2cm] FOURTEST($X_k$) \\
%\tab[2cm] \textbf{if} $3\leq k \leq 5:$\\
%%\tab[3cm] $T = T\cap S$\\
%\tab[3cm] $\text{A}_{k+1}(X_{k+1},T_k\setminus Y)$\\
%\tab \textbf{else}:\\ 
%%\tab[2cm] $T = T\setminus Y$\\
%%\tab[2cm] $X_{k-1} = Y$\\
%\tab[2cm] \textbf{if} $i>3:$\\
%\tab[3cm] $\text{A}_{k-1}(Y,T_k\setminus Y)$ \\
%\tab[2cm] \textbf{else}:\\
%\tab[3cm] $\text{A}_2(Y)$\\
%}}}
%\caption{Procedure \textbf{$\text{A}_k$} for  $3\leq k\leq 6$}
%\label{Ak}
%\end{center} 
%\end{figure}
% 

\begin{algorithm}
	\caption{$\text{A}_2(\text{set } \{x,y,z,w\})$ }
	\label{Single}
				\For{$p\in \{x,y,z\}$}{
				TEST$(\{p\})$  \\
				\eIf{$p$ is defective}{
				$D = D \cup \{p\}$}
				{$ G = G \cup \{p\}$}
				$ S = S \setminus \{p\}$ }
				\If{$\{x,y,z\}$ is good}
				{$S = S \setminus \{w\}$  \\
				$D = D \cup \{w\}$}
				$T_3 = T_3\cap S$\\
				$\text{A}_3(X_3,T_3)$\\

\end{algorithm}

%\begin{figure}
%\begin{center}
%{\footnotesize\colorbox{grey}{\parbox{6.5cm}{
%{\normalsize \textbf{Procedure $\text{A}_2(\text{set } \{x,y,z,w\})$}} \\[0.25cm]
%\textbf{for} $p\in \{x,y,z\}:$\\
%\tab TEST$(\{p\})$  \\
%\tab \textbf{if} $p$ is defective: \\
%\tab[2cm] $D = D \cup \{p\}$  \\
%\tab \textbf{else}:\\
%\tab[2cm] $ G = G \cup \{p\}$ \\
%\tab $ S = S \setminus \{p\}$ \\[0.25cm]
%\textbf{if} $\{x,y,z\}$ is good: \\
%\tab $S = S \setminus \{w\}$  \\
%\tab $D = D \cup \{w\}$  \\
%$T_3 = T_3\cap S$\\
%$\text{A}_3(X_3,T_3)$\\
%}}}
%\caption{Procedure \textbf{$\text{A}_2$}}
%\label{Single}
%\end{center} 
%\end{figure}

Then the algorithm enters an {\it expansion phase} which we explain and analyse first: 
\subsubsection*{Expansion Phase}
With each negative test we expand the size of the test set by a factor of $8$. That means our test sets have sizes 
$64 \cdot 8, 64 \cdot 8^2, 64 \cdot 8^3, \ldots$ until one of the following two cases occurs: Either we get a contaminated test set $X$ of size $64 \cdot 8^k, k \geq 1$ or the remaining set $S$ has less than $64 \cdot 8^k, k \geq 1$ elements. Note that after two negative tests, we test the whole remaining set. If it is good, the algorithm stops. This means that all remaining elements are good
and we needed 3 tests for verifying this. Of course, this can happen only in one expansion phase and the additive constant $a$ takes again care of this situation. 

Assume first that we get a positive test with $|X|= 64 \cdot 8^k, k \geq 1$. By now, we have detected $64 \frac{8^k-1}{7}$ good elements using $k+1$ or $k+2$ tests, depending on whether $k=1$ or not. Now this defective group $X$ runs the subroutine FOURTEST($X$) with $|X|=64\cdot 8^k = 2^l$, $l=3k+6$.
According to the discussion above, the expansion phase ends with  one of the following cases:

\begin{itemize}
\item We detected $d_1=1$ defective element and $g_1=64 \frac{8^k-1}{7}$ good elements using $t_1$ tests, where
$t_1=k+2+1+3k+4 =4k+7$  if $k>1$ and $t_1=4k+6=10$ if $k=1$.
%Note that $1.45198(\log ( 64 \frac{8^k-1}{7}+1) +1.1699)$ is larger than $4k+6$ if $k=1$ and larger than $4k+7$ for $k > 1$. 
\item  We found $g_2:= 64 \frac{8^k-1}{7}+32 \cdot 8^k$ good elements and $d_2=1$ defective ones using at most $t_2=4k+9$ tests
\item $g_3:= 64 \frac{8^k-1}{7}+16 \cdot 8^k$ good elements and $d_3=2$ defective ones using at most
$t_3=7k+13$ tests.

\item We found $g_4:= 64 \frac{8^k-1}{7}+48 \cdot 8^k$  good items and $d_4=1$ defectives by at most $t_4=4k+9$ tests.
\end{itemize}

It is easy to check that 
$$ 1.45198\, d_i \left(\log \frac{g_i+d_i}{d_i}+1.1699 \right)> t_i$$
for $1\leq i\leq 4$. \\
In the second case we have less than $64 \cdot 8^k$ elements left. Since this phase of the algorithm started with more than $64+8\cdot 64$ elements, we have $k\geq 2$ and we know that the remaining set
is defective, hence we can extract a defective element by binary splitting. The number of tests is thus at most $k+2+3k+6=4k+8$ which is smaller than  
$1.45198(\log ( 64 \frac{8^k-1}{7}+1) +1.1699)$ for $k\geq 2$.

\subsubsection*{First Set Contaminated}

Now we consider the case that the first set $X$ of $64$ elements is contaminated. If we would use FOURTEST, it might happen that one defective and no good element are detected by 6 tests which is too much. Instead, we resort to procedure $A_6$ (compare Algorithm \ref{Ak}) which uses $X_6=X$ and $T_6 = S\setminus X$ as a source of new elements. We first test a set $Y$ of 32 new elements from $T_6$. If it is negative,
we use FOURTEST on $X_6$ and proceed with MAIN, hence choosing 64 new elements and so on. If it is positive, procedure $A_5$ (compare Algorithm \ref{Ak}) is called with $X_5=Y.$ In general, for each $k \geq 3$, the procedure $A_k,$ formally described in Algorithm \ref{Ak}, is applied to a defective subset $X_k$ of size $2^k.$\\
The algorithm first tests a set $Y$ of $2^{k-1}$ elements from $T_k$ which is disjoint from $X_k$. If $Y$ is defective, algorithm $A_{k-1}$ is applied to it.

If $Y$ is good, FOURTEST($X_k$) is executed and $A_{k+1}$ is called again if $k<6.$
The result of this execution of $A_k$ is one of the following:
\begin{itemize}
\item $d_1=1$ defective and $g_1=2^{k-1}$ good elements by $t_1= k+1$ tests.
\item $d_2=1$ and $g_2=2^k$ good elements by $t_2=k+3$ tests.
\item $d_3=2$ defectives and $g_3=3\cdot 2^{k-2}$ good elements by $t_3=2k+1$ tests.
\item $d_4=1$ defective and $g_4=5 \cdot 2^{k-2}$ good elements by $t_4=k+3$ tests.
\end{itemize}

$
 t < 1.45198 \, d\left(\log\frac{g+d}{d}+1.1699 \right)
$
which includes the tests of $X_k$ and $Y$.\\
\noindent
It may happen that $A_3$ generates a defective set $Y$ of cardinality $4$ and calls procedure $A_2(Y)$ formally described in Algorithm \ref{Single}.\\
At this point, we test three elements of $Y$ individually and identify at least $3$ elements with $4$ tests. Hence, we are left with the following cases:
\begin{itemize}
\item   one defective and two good elements: \newline				 $1.45198 \cdot 1 \cdot (\log 3+1.1699)>4$ 
										$\Leftrightarrow$ $1.45198>\nolinebreak \frac{4}{\log 3+1.1699} \approx 1.451978$ 
\item   one defective and three good elements:  $1.45198 \cdot 1 \cdot (2+1.1699)>4 $ 
\item   two defectives and one good element:  $1.45198 \cdot  2 \cdot (\log 3-1+1.1699)>5$  
\item   three defectives and no good elements:  $1.45198 \cdot 3 \cdot (\log 1 + 1.1699)> 5$ 
 \end{itemize}

\noindent
The worst over-all case is the first one and so we get $c< 1.45198$ as an upper bound.
Summarizing, we get two estimates for the number of tests used by our algorithm $A$:
\begin{theorem}\label{2}
There is some constant a such that for all $0\leq d <n$ we have
\begin{itemize}
 \item[i)]  $M_A(d\,|\,n) \leq 1.45198  d  \left(\log \frac{n}{d}+1.1699 \right)+a=:f(n,d)$\ ,
 \item[ii)]  $M_A(d\,|\,n)\leq \frac{4}{3}n$\ .
\end{itemize}

\end{theorem}
\noindent
\textit{Proof:}
\begin{enumerate}
\item[i)] The proof is by induction on $d$ and just summarizes the previous discussion. 
If $d=0$, we either need at most $|S|$ individual tests for $|S|\leq 576$ or just 3 tests. 

\noindent
So assume that $d>0$. If the algorithm starts with a negative test, it enters the expansion phase
and identifies $g'$ good and $d'>0$ defective elements by $t$ tests with 
$$t <1.45198 d' \left(\log \frac{g'+d'}{d'}+ 1.1699 \right)$$ 
as described above. By induction, the total number of tests is at most  
$$1.45198 d' \left(\log\frac{g'+d'}{d'}+ 1.1699 \right)+ f(n-(g'+d')),d-d')\leq f(n,d)$$
by Theorem \ref{Jensen}.

Otherwise, the algorithm calls procedure $A_6$. By the discussion above and similar reasoning, inequality \textit{i)} holds in this case as well.

\item[ii)] This follows by induction on $n$ and is trivial for small $n$ or if the algorithm first enters the expansion phase. If algorithm $A_k$, $k>2$, chooses some good set $Y$ with
$|Y|=2^{k-1}$, then FOURTEST($X$) is executed. It is easy to check that the number of tests of procedure FOURTEST plus the tests for $X$ and $Y$ are less than $|Y|$ plus the
number of elements identified by FOURTEST. Finally, if $A_2$ is executed, at least three elements are identified by four tests, the case which is responsible for the ratio $4/3$ in inequality
\textit{ii)}.

\end{enumerate}\hfill $\Box$

\noindent
This leads to:
\begin{theorem}\label{mainresult}
Algorithm $A$ is $1.452$-competitive.
\end{theorem}

\noindent
\textit{Proof:}\\
We may assume that $d>0$ and $n$ is sufficiently large. Moreover, we consider three different ranges of $r=n/d.$
If $r<2,$ we use Theorem \ref{2} ii) and Theorem \ref{Corollary:Cutoff}:
$$\frac{M_A(d\,|\,n)}{M(d,n)}\leq \frac{\frac{4}{3}n}{n-1} \longrightarrow \frac{4}{3} < 1.452 \ .$$
For $3 \geq r \geq 2,$ we apply Theorem \ref{2} ii) and Theorem \ref{Theorem: Stirling}:
$$\frac{M_A(d\,|\,n)}{M(d,n)}\leq \frac{\frac{4}{3}n}{\frac{n}{3}\left(\log 3+2 (\log 3-\log 2) \right)+O(\log n)} \longrightarrow \frac{4}{3\log 3 - 2} < 1.452 \ .$$
Finally, if $r \geq 3$, we have with Theorem \ref{2} i) and Theorem \ref{Theorem: Stirling}:
$$\frac{M_A(d\,|\,n)}{M(d,n)}\leq \frac{1.45198  d \left(\log \frac{n}{d} +1.1699\right)+a}
{d  \left( \log \frac{n}{d}+1.1699\right)- \frac{1}{2}\log d-\frac{3}{2}} \longrightarrow 1.45198<1.452 \ . $$ 
\hfill $\Box$

%%%%%%%%%%%%%%%%%%%%%%%%%%%%%%%%%%%%%%%%%%%%%%%%%%%%%%%%%%%%%%%%%%%%%%%%%%%%%%%%%%%%%%%%%%%%%%%%%%%%%%%%%%%%%%%%%%%%%%%%%%%%%%%%
\section{A lower bound}

To prove a lower bound for competitive constants $c$, we assume familiarity with the definition of the worst case search length as a game between two players
as developed, e.g., in \cite[section~1.10]{Aigner1988}: %(Aigner, section 1.10)
 The first player (Algy) chooses a test set $W$ and the second (the Strategist, his adversary) answers \textit{good} or \textit{defective}.
Of course, his answers have to be consistent and his goal is to force Algy to ask as many questions as possible whereas Algy
wants to minimize that number. Moreover, we can assume that Algy's test sets do not contain any classified elements.
The worst case number of tests is just the number of tests in the game if both players play optimally from their point of view.
In particular, if the adversary player follows some specific strategy and Algy plays optimally against this strategy, the number of tests needed is a lower bound
for the general worst case number of tests. We are going to define and analyze such a strategy ${\cal S}$.

\noindent
In order to precisely define the answers of the Strategist in the $t$-th round, we have to keep track of the results of the preceding tests.
This is done by some hypergraph $H_t=(S_t, E_t)$, where $S_t$ is the set of all items which have not been identified as good or defective after $t$ tests, and $E_t$ consists of all inclusion-minimal
sets which are known to contain a defective element after $t$ tests. $H_t$ is defined inductively as follows:
\begin{itemize}
\item $H_0:=(S,\emptyset)$
\item If the hypergraph $H_{t-1}$ is defined (after $t-1$ tests and answers), assume that Algy chooses some test set $W,$ w.l.o.g. $\emptyset\neq W\subset S_{t-1}$.
\begin{itemize}
\item If the adversary answers \textit{good}, all elements of $W$ are removed from $S_{t-1}$. Furthermore, let $U$  denote the set of items $x\in S_{t-1}$
such that $\{x\}=Y\setminus W$ for some $Y\in E_{t-1},$ i.e., all vertices which lie in edges of size one after the removal of $W.$ The elements of $U$ are identified as \textit{defective} and are removed from $S_{t-1}$ as well. Thus
$S_t:=S_{t-1}\setminus (U\cup W)$. Now let $E_t':=\{Y\setminus W| Y\in E_{t-1} \text{ and }Y\cap U=\emptyset\}$. $E_t$ is defined as the set of all inclusion-minimal sets in $E_t'$.

\item If the adversary answers \textit{defective}, we distinguish two cases:
If $|W|>1$, we let $S_t:=S_{t-1}$ and $E_t$ consists of all inclusion-minmal sets in $E_{t-1}\cup \{W\}$. 
Otherwise, if $|W|=1,$ set $S_t:=S_{t-1}\setminus W$ and $E_t':=\{Y| Y\in E_{t-1} \text{ and }Y\cap W=\emptyset\}$. $E_t$ is defined as the set of all inclusion-minimal sets in $E_t'$.
\end{itemize}
\end{itemize}
Recall that a \textit{vertex cover} of $H_t$ is a set $X$ with $X\cap Y\neq \emptyset$ for all  $Y\in E_t$.
After $t$ tests, the possible sets of defectives are just the unions of the defective elements in $S\setminus S_t$ with the vertex covers of $H_t$ and with elements which are not covered by any edge in $H_t$. For $|S_t|\geq 2$ all sets in $E_t$ have cardinality at least two. Hence, subsets of $S_t$ with $|S_t|-1$ elements are vertex covers. Therefore, the set of defectives is not determined yet.
The case $|S_t|=1$ arises only if $E_t=\emptyset$. It follows that the set of defectives is determined iff $S_t=\emptyset$.

\noindent
We can now define our strategy ${\cal S}$:

\begin{itemize}
\item[-] If a set $W$ of at least two elements is tested, the answer is ``defective".
\item[-] If a singleton is tested which is not contained in some $Y\in E_t$, the answer is ``good".
\item[-] If a singleton is tested which is contained in  some $Y\in E_t$, the answer is ``defective".
\end{itemize}
\noindent

\begin{theorem}\label{strategy}
If answers are given according to strategy ${\cal S}$, the following assertions hold for each group testing algorithm $A$:
\begin{itemize}

\item[i)] All elements will never be declared defective.
\item[ii)] If the number of defectives in the game is $d$, then $A$ used at least $n+d$ tests.
\item[iii)] If $A$ is $c$-competitive, we have for $n\to \infty$
\begin{align*}
 c \geq \frac{M_A(d\,|\, n)}{M(d,n)} -o(1)\geq \min\left\lbrace\left. \frac{n+d}{M(d,n)}\right|0 \leq d < n\right\rbrace -o(1).
\end{align*}
\end{itemize}
\end{theorem}

\noindent
\textit{Proof:}\\
$i)$ Assume to the contrary that all elements are declared defective. Then the Strategist always answers \textit{defective}. As long as $|S_t|\geq 2$, a new defective can only be identified if
Algy chooses a singleton as test set. Hence, for some $t$, we have $|S_t|=1$ and $E_t=\emptyset$. But then the next test has to be $W=S_t$ and the answer is \textit{good},
a contradiction.

\noindent
$ii)$ It is clear that the $(n-d)$ good items are identified by testing $(n-d)$ singletons. If some item $x$ is declared defective in the $(t+1)$-th test, the corresponding test has to be $W=\{x\}$ and there has to be some prior test of a set
$Y=Y(x)$ such that $x\in Y \in E_t$. Since $Y\notin E_{t+1}$, $Y(x)\neq Y(x')$ for all other defective elements $x'\neq x$. It follows that $A$ needs at least two tests for each defective element.
The total number of tests is thus at least $2d+(n-d)=n+d$.

\noindent
$iii)$ Obvious, by $i)$ and $ii)$.\hfill $\Box$

In the range $d \geq \frac{16}{43}n$, individual tests are optimal by Theorem \ref{Leu2002}, hence 
\begin{align*}
 \frac{n+d}{M(d,n)} = \frac{n+d}{n-1} \geq \frac{59}{43}-o(1)> 1.37
\end{align*}

For $d<\frac{16}{43}n$, we estimate
\noindent 
\begin{align*}
\frac{n+d}{M(d,n)}\geq \frac{n+d}{M_{A_d}(d,n)} \ .
\end{align*}
where $A_d$ is Allemann's group testing algorithm from Theorem \ref{Allemann2013}, hence
\begin{align*}
 M_{A_d}(d,n)\leq \log \left(\begin{pmatrix}
       n \\d
      \end{pmatrix}\right)
+0.255d +\frac{1}{2} \log (d)+6.5
\end{align*}

\begin{figure}[t]
\begin{center}
 \includegraphics[width= 7cm]{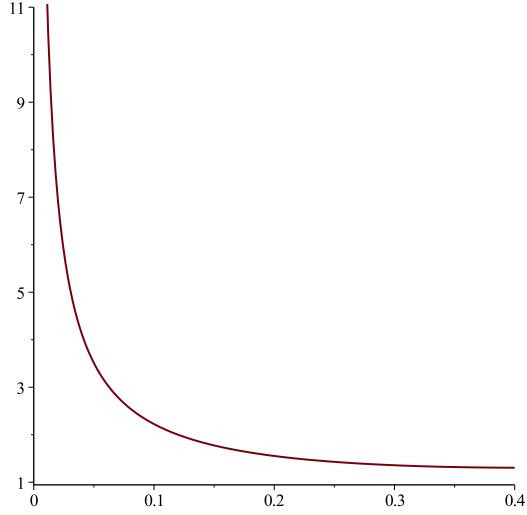}
\end{center}
\caption{Function $f(p)$}
\label{Skizze}
\end{figure}

\noindent
By using the estimate 
$${n\choose d}\leq \Big(e\frac{n}{d}\Big)^d,$$
it is easy to see that for each constant B, there exists some $\epsilon>0$ such that 
$$\liminf_{n\to \infty,\, d\leq \epsilon n}\frac{n+d}{M_{A_d}(d,n)}\geq B:$$
 
Now let $p:=p_n:= \frac{d}{n}$
and $h(p)$ $= -p \log p-(1-p) \log (1-p)$ 
 denote the binary entropy function. It is well known that $\log \left({n \choose d}\right)
= n(h(p))+o(1))$ for $\epsilon\leq p\leq 1/2$. We may thus write $n(h(p))+0.255p+o(1)$ for the number of tests in algorithm $A_d$, hence
\begin{align*}
 \frac{n+d}{n(h(p))+0.255d+o(1)}= \frac{1+p}{h(p)+0.255p+o(1)} \ .
\end{align*}

\noindent
The function $f(p):=\frac{1+p}{h(p)+0.255p}$ is non-increasing on $(0,0.4)$ (see Figure \ref{Skizze}). 
For $p= \frac{16}{43}$ we get $f(p)\approx 1.3103$.
This shows our final result:

\begin{theorem}
 If an algorithm is c-competitive for sufficiently large $n$, then we have $c > 1.31$.
\end{theorem}

\section{Conclusion}
We do not have a serious conjecture for the best possible competitive constant ($4/3$ would be nice) and guess that both bounds we developed can be improved significantly.

 \bibliographystyle{unsrt}
\bibliography{QuellenBounds}

\end{document}